\documentclass[12pt,english]{article}
\usepackage{babel}
\usepackage[cp1251]{inputenc}
\usepackage{amsmath,amssymb}
\usepackage{bm}

\begin{document}

\newpage

\begin{center}
{\bf \Large
Estimates for solutions in a predator--prey model
describing plankton--fish interaction\footnote{The study was carried out
within the framework of the state contract of the Sobolev Institute of Mathematics 
(project no.~FWNF-2022-0008).}}
\end{center}

\begin{center}
{\bf M.A.~Skvortsova$^{1,2}$}
\end{center}

\begin{center}
${}^1$Sobolev Institute of Mathematics, \\
${}^2$Novosibirsk State University, \\
Novosibirsk, Russia
\end{center}

\begin{center}
sm-18-nsu@yandex.ru
\end{center}

{\bf Abstract.}
In the paper we consider a system of differential equations
with two delays describing plankton--fish interaction.
We study stability of the equilibrium point corresponding
to the presence of phytoplankton and zooplankton in the system
and the absence of fish.
In the case of asymptotic stability,
we indicate conditions for the initial data
under which solutions stabilize to this equilibrium point
and establish estimates for solutions characterizing
the stabilization rate at infinity.
The results are obtained using Lyapunov--Krasovskii functional.
The obtained theoretical results can be used for numerical study
of behavior of solutions to the considered model.

\vskip10pt

\linebreak
equations, equilibrium point,
asymptotic stability, estimates for solutions,
attraction set, Lyapunov--Krasovskii functional

\vskip10pt



\begin{center}
{\bf 1. Introduction}
\end{center}

At present, there exist a large number of works devoted to the study
of biological models described by delay differential equations
(see, for example, the monographs~\cite{G1992, K1993, S1995, E2009}
and bibliography therein).
Among studied models, the models of population dynamics are widespread,
in particular, models of predator--prey type, based on the classical
predator--prey model that was introduced independently
by A.J.~Lotka~\cite{L1925} and V.~Volterra~\cite{V1926}.
An overview of some results for predator--prey models
with delay is contained, for example, in~\cite{LChA2002, R2009}.

In particular, predator--prey models are used
when describing plankton--fish interaction
(see, for example,~\cite{PCh2015, KCP2016, MW2018, ZhS2019, RTM2020, ThO2020, S2021}).
Taking into account a model proposed in~\cite{ThO2020},
in the present paper, we study the model of the following form:
$$
\left\{
\begin{array}{l}
\displaystyle
\frac{d}{dt}x(t)=rx(t) \left( 1-\frac{x(t)}{K} \right)-c_1 x(t) y(t),
\\
\\
\displaystyle
\frac{d}{dt}y(t)=-d_1 y(t)+e_1 c_1 x(t-\tau_1) y(t-\tau_1)-c_2 y(t) z(t),
\\
\\
\displaystyle
\frac{d}{dt}z(t)=-d_2 z(t)+e_2 c_2 y(t-\tau_2) z(t-\tau_2),
\end{array}
\right.
\eqno (1)
$$
which can be also considered as a model of plankton--fish interaction
(note that, in~\cite{ThO2020}, it was considered the case
$\tau_2=0$,
but nonlinear terms had a more general form).
In this system,
$x(t)$
is the amount of phytoplankton,
$y(t)$
is the amount of zooplankton, and
$z(t)$
is the number of fish.
It is assumed that phytoplankton is the favorite food of zooplankton,
which serves as the favorite food of fish.
The delay parameter
$\tau_1 \ge 0$
is responsible for time required for the appearance of new zooplankton,
and the delay parameter
$\tau_2 \ge 0$
is responsible for time of fish maturation.
The coefficients of the system have the following meaning:
$r>0$ is intrinsic growth rate of phytoplankton,
$K>0$ is environmental carrying capacity of phytoplankton,
$d_1>0$ is mortality rate of zooplankton,
$d_2>0$ is mortality rate of fish,
$c_1 \ge 0$ is predation rate of zooplankton,
$c_2 \ge 0$ is predation rate of fish,
$e_1=b_1 e^{-c_1 \tau_1}$,
$b_1 \ge 0$ is birth rate of zooplankton,
$e_2= b_2 e^{-c_2 \tau_2}$,
$b_2 \ge 0$ is birth rate of fish.

We consider system~(1) for
$t>0$,
assuming that the initial conditions are given on the segment
$\theta \in [-\tau_{\max},0]$,
$\tau_{\max}=\max\{ \tau_1, \tau_2 \}$:
$$
\left\{
\begin{array}{l}
\displaystyle
x(\theta)=\varphi(\theta) \ge 0, \quad \theta \in [-\tau_1,0], \quad x(+0)=\varphi(0)>0,
\\
\\
\displaystyle
y(\theta)=\psi(\theta) \ge 0, \quad \theta \in [-\tau_{\max},0], \quad y(+0)=\psi(0),
\\
\\
\displaystyle
z(\theta)=\eta(\theta) \ge 0, \quad \theta \in [-\tau_2,0], \quad z(+0)=\eta(0),
\end{array}
\right.
\eqno (2)
$$
where
$\varphi(\theta)$,
$\psi(\theta)$,
$\eta(\theta)$
are continuous functions.
It is well known that a solution to the initial value problem~(1), (2)
exists and is unique.
Moreover, by analogy with~\cite{ThO2020}, it is not difficult to show
that the solution is defined on the entire right half-axis
$\{ t>0 \}$,
has non-negative components,
and each component of the solution is a bounded function.
In other words, there exist constants
$M_1, \ M_2, \ M_3>0$
such that, for all
$t>0$,
the inequalities are valid
$$
0 \leq x(t) \leq M_1, \quad
0 \leq y(t) \leq M_2, \quad
0 \leq z(t) \leq M_3,
$$
i.e., the amount of plankton and fish cannot increase indefinitely.

We write out all equilibrium points of system~(1)
with non-negative components.

1) If
$$
d_1 \ge e_1 c_1 K,
$$
then there are two equilibrium points in the system:
$(x(t),y(t),z(t))=(0,0,0)$
and
$(x(t),y(t),z(t))=(K,0,0)$.
Equilibrium point
$(0,0,0)$
corresponds to the complete extin\-ction of fish and plankton;
equilibrium point
$(K,0,0)$
corresponds to the presence of only phytoplankton in the system
and the absence of zooplankton and fish.

2) If
$$
e_1 c_1 K \left( 1-\frac{c_1 d_2}{e_2 c_2 r} \right) \leq d_1<e_1 c_1 K,
$$
then there are three equilibrium points in the system:
$(0,0,0)$,
$(K,0,0)$,
and
$(x_0,y_0,0)$,
where
$$
x_0=\frac{d_1}{e_1 c_1}, \quad
y_0=\frac{r}{c_1} \left( 1-\frac{d_1}{e_1 c_1 K} \right).
\eqno (3)
$$
Equilibrium point
$(x_0,y_0,0)$
corresponds to the presence of phytoplankton and zooplankton
in the system and the absence of fish.

3) If
$$
d_1<e_1 c_1 K \left( 1-\frac{c_1 d_2}{e_2 c_2 r} \right),
$$
then there are four equilibrium points in the system:
$(0,0,0)$,
$(K,0,0)$,
$(x_0,y_0,0)$,
and
$(x^{*},y^{*},z^{*})$,
where
$$
x^{*}=K \left( 1-\frac{c_1 d_2}{e_2 c_2 r} \right),
\quad
y^{*}=\frac{d_2}{e_2 c_2},
\quad
z^{*}=\frac{1}{c_2}
\left( e_1 c_1 K \left( 1-\frac{c_1 d_2}{e_2 c_2 r} \right)-d_1 \right).
$$
Equilibrium point
$(x^{*},y^{*},z^{*})$
corresponds to the coexistence
of fish, phytoplankton and zooplankton in the system.

The question of stability of equilibrium points was considered in~\cite{ThO2020} when
$\tau_2=0$.
We formulate the corresponding results.

1) For
$\tau_2=0$,
equilibrium point
$(0,0,0)$
is unstable~\cite{ThO2020}.
For
$\tau_2>0$,
by analogy with~\cite{ThO2020}, it is also not difficult to show
that equilibrium point
$(0,0,0)$
will be unstable.

2) For
$\tau_2=0$,
if
$d_1>e_1 c_1 K$,
then equilibrium point
$(K,0,0)$
is asymptotically stable, and if
$d_1<e_1 c_1 K$,
then equilibrium point
$(K,0,0)$
is unstable~\cite{ThO2020}.
For
$\tau_2>0$,
the result also remains valid:
if
$d_1>e_1 c_1 K$,
then
$(K,0,0)$
is asymptotically stable, and if
$d_1<e_1 c_1 K$,
then
$(K,0,0)$
is unstable.

3) For
$\tau_1=\tau_2=0$,
if
$
\displaystyle
e_1 c_1 K \left( 1-\frac{c_1 d_2}{e_2 c_2 r} \right)<d_1<e_1 c_1 K
$,
then equilibrium point
$(x_0,y_0,0)$
is asymptotically stable, and if
$
\displaystyle
d_1<e_1 c_1 K \left( 1-\frac{c_1 d_2}{e_2 c_2 r} \right)
$,
then equilibrium point
$(x_0,y_0,0)$
is unstable~\cite{ThO2020}.

4) For
$\tau_1=\tau_2=0$,
if
$
\displaystyle
d_1<e_1 c_1 K \left( 1-\frac{c_1 d_2}{e_2 c_2 r} \right)
$,
then equilibrium point
$(x^{*},y^{*},z^{*})$
is asymptotically stable. For
$\tau_1>0$
and
$\tau_2=0$,
conditions of the asymptotic stability of equilibrium point
$(x^{*},y^{*},z^{*})$
depend on the delay parameter
$\tau_1$~\cite{ThO2020}.

Note that along with the study of stability of equilibrium points,
an important issue is also finding acceptable conditions for initial data
under which solutions are stabilized, and obtaining estimates
characterizing the stabilization rate of solutions at infinity.
When studying these questions, Lyapunov--Krasovskii functionals
are actively used in the litera\-ture
(see, for example, \cite{DM2005, KhIK2005, MKh2005, DM2007, D2009_1,
DM2014, DM2015_1, DM2015_2, M2017_1, M2017_2, DMS2019,
M2020_1, M2020_2, Y2020, M2021}
for general classes of systems with delay,
\cite{S2021, S2016, S2018_1, S2018_2, S2018_3, S2019, SY2021}
for certain biological models).

For equilibrium point
$(K,0,0)$
of system~(1), the question of obtaining estimates
characterizing the stabilization rate of solutions
to this equilibrium point
was considered in~\cite{S2021}.
In the work~\cite{S2021}, under the condition
$d_1>e_1 c_1 K$,
which guarantees the asymptotic stability of equilibrium point
$(K,0,0)$,
it was shown that all solutions to system~(1)
with initial data of the form~(2) are stabilized
at infinity to this equilibrium point.
At the same time, estimates were established
for all components of the solution to the initial value problem~(1), (2),
characterizing the stabilization rate of the solution at infinity to equilibrium point
$(K,0,0)$.

The aim of this work is to study the stability of equilibrium point
$(x_0,y_0,0)$
for
$\tau_1 \ge 0$
and
$\tau_2\ge 0$,
finding conditions on the initial data under which solutions are stabilized
to this equilibrium point, and obtaining estimates of solutions
characterizing the stabilization rate at infinity.

For equilibrium point
$(x^{*},y^{*},z^{*})$,
the issue of global asymptotic stability
was considered in~\cite{ThO2020} when
$\tau_2=0$.
The question of obtaining estimates of the stabilization rate
of solutions to this equilibrium point is still open.


\begin{center}
{\bf 2. Stability of equilibrium point $(x_0,y_0,0)$}
\end{center}

In this section, we obtain the stability conditions for equilibrium point
$(x_0,y_0,0)$,
where
$x_0$
and
$y_0$
are defined in~(3).
We assume that the inequality is satisfied
$$
d_1 \leq e_1 c_1 K,
$$
which guarantees the non-negativity of components
$x_0$
and
$y_0$.
First, we reduce the problem of stability of equilibrium point
$(x_0,y_0,0)$
of system~(1) to study the stability of the zero solution.
To do this, we change the variables
$$
x(t)=x_0+\widetilde{x}(t),
\quad
y(t)=y_0+\widetilde{y}(t),
\quad
z(t)=\widetilde{z}(t).
\eqno (4)
$$
Then the system will take the form
$$
\frac{d}{dt} \widetilde{\bf y}(t)
=A\widetilde{\bf y}(t)
+B_1 \widetilde{\bf y}(t-\tau_1)
+B_2 \widetilde{\bf y}(t-\tau_2)
$$
$$
+F(\widetilde{\bf y}(t))
+G_1(\widetilde{\bf y}(t-\tau_1))
+G_2(\widetilde{\bf y}(t-\tau_2)),
\eqno (5)
$$
where
$$
\widetilde{\bf y}(t)
=\begin{pmatrix}
\widetilde{x}(t) \\ \widetilde{y}(t) \\ \widetilde{z}(t)
\end{pmatrix},
\quad
A=\begin{pmatrix}
-\frac{r}{K} x_0 & -\frac{d_1}{e_1} & 0 \\
0 & -d_1 & -c_2 y_0 \\
0 & 0 & -d_2
\end{pmatrix},
\eqno (6)
$$
$$
B_1=\begin{pmatrix}
0 & 0 & 0 \\
e_1 c_1 y_0 & d_1 & 0 \\
0 & 0 & 0
\end{pmatrix},
\quad
B_2=\begin{pmatrix}
0 & 0 & 0 \\
0 & 0 & 0 \\
0 & 0 & e_2 c_2 y_0
\end{pmatrix},
\eqno (7)
$$
$$
F(\widetilde{\bf y}(t))
=\begin{pmatrix}
-\frac{r}{K} \widetilde{x}^2(t)-c_1 \widetilde{x}(t) \widetilde{y}(t) \\
-c_2 \widetilde{y}(t) \widetilde{z}(t) \\
0
\end{pmatrix},
\eqno (8)
$$
$$
G_1(\widetilde{\bf y}(t-\tau_1))
=\begin{pmatrix}
0 \\
e_1 c_1 \widetilde{x}(t-\tau_1) \widetilde{y}(t-\tau_1) \\
0
\end{pmatrix},
\eqno (9)
$$
$$
G_2(\widetilde{\bf y}(t-\tau_2))
=\begin{pmatrix}
0 \\
0 \\
e_2 c_2 \widetilde{y}(t-\tau_2) \widetilde{z}(t-\tau_2)
\end{pmatrix}.
\eqno (10)
$$

To obtain stability conditions for the zero solution to system~(5),
we use the theorem of stability in the first approximation
(see, for example,~\cite{K1959}, chapter~7, section~33):

1) if all roots of characteristic quasi-polynomial
$$
Q(\lambda)=\det(\lambda E-A-e^{-\lambda\tau_1}B_1-e^{-\lambda\tau_2}B_2)=0
\eqno (11)
$$
are contained in the left half-plane
$\mathbb{C}_{-}=\{ \lambda \in \mathbb{C}: {\rm Re}\,\lambda<0 \}$,
then the zero solution to system~(5) is asymptotically stable,
so equilibrium point
$(x_0,y_0,0)$
of system~(1) is also asymptotically stable;

2) if there is a root of equation~(11) lying in the right half-plane
$\mathbb{C}_{+}=\{ \lambda \in \mathbb{C}: {\rm Re}\,\lambda>0 \}$,
then the zero solution to system~(5) is unstable,
hence equilibrium point
$(x_0,y_0,0)$
is unstable too.

We formulate the corresponding result.

{\bf Lemma.}
{\it
1) Under the condition
$$
e_1 c_1 K \cdot \max\left\{
\frac{1}{3}, \left( 1-\frac{c_1 d_2}{e_2 c_2 r} \right)
\right\}<d_1<e_1 c_1 K,
\eqno (12)
$$
equilibrium point
$(x_0,y_0,0)$
of system~$(1)$ is asymptotically stable.

2) If the inequality holds
$$
d_1<e_1 c_1 K \left( 1-\frac{c_1 d_2}{e_2 c_2 r} \right),
\eqno (13)
$$
then equilibrium point
$(x_0,y_0,0)$
of system~$(1)$ is unstable.
}

{\bf Remark.}
Under condition~(12), there are three equilibrium points
of the system~(1) with non-negative components;
under condition~(13), there are four equilibrium positions.

{\bf Proof.}
We write down the characteristic quasi-polynomial:
$$
Q(\lambda)=\det(\lambda E-A-e^{-\lambda\tau_1}B_1-e^{-\lambda\tau_2}B_2)
$$
$$
=\left|
\begin{matrix}
\lambda+\frac{r}{K} x_0 & \frac{d_1}{e_1} & 0 \\
-e_1 c_1 y_0 e^{-\lambda\tau_1} & \lambda+d_1-d_1 e^{-\lambda\tau_1} & c_2 y_0 \\
0 & 0 & \lambda+d_2-e_2 c_2 y_0 e^{-\lambda\tau_2}
\end{matrix}
\right|
$$
$$
=Q_1(\lambda) Q_2(\lambda)=0,
$$
where
$$
Q_1(\lambda)=\left( \lambda+\frac{r}{K} x_0 \right)
(\lambda+d_1-d_1 e^{-\lambda\tau_1})
+c_1 d_1 y_0 e^{-\lambda\tau_1},
$$
$$
Q_2(\lambda)=\lambda+d_2-e_2 c_2 y_0 e^{-\lambda\tau_2}.
$$

First, we consider quasi-polynomial
$Q_2(\lambda)$.
It is well known (see, for example,~\cite{EN1971}, chapter~3, section~3)
that under the condition
$0 \leq e_2 c_2 y_0<d_2$,
all roots of quasi-polynomial
$Q_2(\lambda)$
are contained in the left half-plane
$\mathbb{C}_{-}$,
and under the condition
$d_2<e_2 c_2 y_0$,
there is a root of quasi-polynomial
$Q_2(\lambda)$
belonging to the right half-plane
$\mathbb{C}_{+}$.
Taking into account formula~(3) for the value
$y_0$,
we obtain the following result:

1) if
$$
e_1 c_1 K \left( 1-\frac{c_1 d_2}{e_2 c_2 r} \right)
<d_1 \leq e_1 c_1 K,
\eqno (14)
$$
then all roots of quasi-polynomial
$Q_2(\lambda)$
are contained in the left half-plane
$\mathbb{C}_{-}$,

2) if condition~(13) is met:
$$
d_1<e_1 c_1 K \left( 1-\frac{c_1 d_2}{e_2 c_2 r} \right),
$$
then there is a root of quasi-polynomial
$Q_2(\lambda)$
belonging to the right half-plane
$\mathbb{C}_{+}$.

Now we consider quasi-polynomial
$Q_1(\lambda)$.
Because
$\lambda=-\frac{r}{K} x_0$
is not the root of this quasi-polynomial,
then it can be rewritten in the form
$$
Q_1(\lambda)=\left( \lambda+\frac{r}{K} x_0 \right)
\left[
\lambda+d_1-d_1 e^{-\lambda\tau_1}
\left( 1-\frac{c_1 y_0}{\left( \lambda+\frac{r}{K} x_0 \right)} \right)
\right]
=\left( \lambda+\frac{r}{K} x_0 \right) Q_3(\lambda).
$$
We find the conditions under which all roots of quasi-polynomial
$Q_3(\lambda)$
are contained in the left half-plane
$\mathbb{C}_{-}$.

Suppose that there is a root
$\lambda^{*} \in \mathbb{C}$
of quasi-polynomial
$Q_3(\lambda)$
such that
${\rm Re}\,\lambda^{*} \ge 0$
and the inequality is fulfilled
$$
\left| 1-\frac{c_1 y_0}{\left( \lambda^{*}+\frac{r}{K} x_0 \right)} \right|<1.
\eqno (15)
$$
Then
$$
|Q_3(\lambda^{*})|
=\left|
\lambda^{*}+d_1-d_1 e^{-\lambda^{*}\tau_1}
\left( 1-\frac{c_1 y_0}{\left( \lambda^{*}+\frac{r}{K} x_0 \right)} \right)
\right|
$$
$$
\ge |\lambda^{*}+d_1|-d_1 e^{-\tau_1 \, {\rm Re}\,\lambda^{*}}
\left| 1-\frac{c_1 y_0}{\left( \lambda^{*}+\frac{r}{K} x_0 \right)} \right|
$$
$$
>|\lambda^{*}+d_1|-d_1
=\sqrt{\left( {\rm Re}\,\lambda^{*}+d_1 \right)^2
+\left( {\rm Im}\,\lambda^{*} \right)^2}-d_1 \ge 0,
$$
i.e.,
$|Q_3(\lambda^{*})|>0$.
This contradicts the fact that
$\lambda^{*}$
is a root of quasi-polynomial
$Q_3(\lambda)$.
Thus, condition~(15) is a sufficient condition
that all roots of quasi-polynomial
$Q_3(\lambda)$
are contained in the left half-plane
$\mathbb{C}_{-}$.

We transform condition~(15):
$$
\left| \lambda^{*}+\frac{r}{K} x_0-c_1 y_0 \right|^2
<\left| \lambda^{*}+\frac{r}{K} x_0 \right|^2,
$$
$$
\left( {\rm Re}\, \lambda^{*}+\frac{r}{K} x_0-c_1 y_0 \right)^2
+({\rm Im}\, \lambda^{*})^2
<\left( {\rm Re}\, \lambda^{*}+\frac{r}{K} x_0 \right)^2
+({\rm Im}\, \lambda^{*})^2,
$$
$$
(c_1 y_0)^2
-2c_1 y_0 \left( {\rm Re}\, \lambda^{*}+\frac{r}{K} x_0 \right)
<0,
$$
$$
0<c_1 y_0<2{\rm Re}\, \lambda^{*}+\frac{2r}{K} x_0.
$$
It follows that when the inequalities are met
$$
0<c_1 y_0<\frac{2r}{K} x_0,
$$
all roots of quasi-polynomial
$Q_3(\lambda)$
are contained in the left half-plane
$\mathbb{C}_{-}$.
By virtue of definitions~(3) of values
$x_0$
and
$y_0$,
this condition is equivalent to the following:
$$
\frac{1}{3} e_1 c_1 K<d_1<e_1 c_1 K.
\eqno (16)
$$

Considering all of the above, we get the final result.

If conditions~(14) and~(16) are fulfilled simultaneously,
then all roots of characteristic quasi-polynomial
$Q(\lambda)$
are contained in the left half-plane
$\mathbb{C}_{-}$.
According to the theorem of stability in the first approximation,
it follows the asymptotic stability of equilibrium point
$(x_0,y_0,0)$
of system~(1).

When condition~(13) is met, characteristic quasi-polynomial
$Q(\lambda)$
has roots lying in the right half-plane
$\mathbb{C}_{+}$.
In this case, equilibrium point
$(x_0,y_0,0)$
of system~(1) is unstable.

Lemma is proved.

{\bf Remark.}
Note that in the case
$$
e_1 c_1 K \left( 1-\frac{c_1 d_2}{e_2 c_2 r} \right) \leq d_1
\leq e_1 c_1 K \cdot \max\left\{
\frac{1}{3}, \left( 1-\frac{c_1 d_2}{e_2 c_2 r} \right)
\right\},
$$
conditions of stability of equilibrium point
$(x_0,y_0,0)$
depend on the delay parameter
$\tau_1$.
This case is not considered in this article.

Throughout the following, we will assume that condition~(12) is satisfied,
which guaran\-tees the asymptotic stability of the considered equilibrium point
$(x_0,y_0,0)$.


\begin{center}
{\bf 3. Construction of the Lyapunov--Krasovskii functional}
\end{center}

The aim of this section is to construct
a Lyapunov--Krasovskii functional.
Hereinafter, with the help of this functional,
conditions for the initial data will be obtained,
under which solutions to system~(1) will stabilize to equilibrium point
$(x_0,y_0,0)$,
and estimates characterizing the stabilization rate at infinity will be established.

When constructing the functional, we use the idea proposed in~\cite{DM2005}.
First, we give the results from this work.
Consider the system of linear delay differential equations
$$
\frac{d}{dt}{\bf y}(t)=A{\bf y}(t)+B{\bf y}(t-\tau), \quad t>0.
\eqno (17)
$$
Suppose that there exist Hermitian positive definite matrices
$H$
and
$K(s) \in C^1([0,\tau])$,
i.e.,
$H=H^{*}>0$
and
$K(s)=K^{*}(s)>0$
are such that
$\frac{d}{ds}K(s)<0$,
$s \in [0,\tau]$,
and matrix
$$
C=-\begin{pmatrix}
HA+A^{*}H+K(0) && HB \\
B^{*}H && -K(\tau)
\end{pmatrix}
$$
is positive definite.
It was shown in~\cite{DM2005} that under these conditions,
the zero solution to system~(17) is asymptotically stable.

The proof of this result is based on the use
of the Lyapunov--Krasovskii functional~\cite{DM2005}
$$
V(t,{\bf y})=\left\langle H{\bf y}(t),{\bf y}(t) \right\rangle
+\int\limits_{t-\tau}^{t}
\left\langle K(t-s){\bf y}(s),{\bf y}(s) \right\rangle ds,
$$
the derivative of which by virtue of system~(17) has the form
$$
\frac{d}{dt}V(t,{\bf y})
=-\left\langle
C
\begin{pmatrix} {\bf y}(t) \\ {\bf y}(t-\tau) \end{pmatrix},
\begin{pmatrix} {\bf y}(t) \\ {\bf y}(t-\tau) \end{pmatrix}
\right\rangle
+\int\limits_{t-\tau}^{t}
\left\langle \frac{d}{dt}K(t-s){\bf y}(s),{\bf y}(s) \right\rangle ds.
$$
It is important to note that in~\cite{DM2005} with the help of this functional
in addition to conditions for the asymptotic stability of the zero solution,
estimates of solutions to system~(17) characterizing the decay rate at infinity were specified.
In this paper, along with the linear system~(17), a nonlinear system was considered,
for which sufficient conditions for the asymptotic stability of the zero solution were also obtained,
conditions for the initial data under which solutions decrease were specified,
and estimates of the rate of solution decrease at infinity were established.

By analogy with the results of the work~\cite{DM2005},
it is possible to obtain sufficient conditions for the asymptotic stability
of the zero solution to linear system of differential equations with two delays
$$
\frac{d}{dt}{\bf y}(t)=A{\bf y}(t)
+B_1{\bf y}(t-\tau_1)+B_2{\bf y}(t-\tau_2), \quad t>0.
\eqno (18)
$$
The following statement is valid.
If there exist matrices
$H=H^{*}>0$,
$K_1(s) \in C^1([0,\tau_1])$,
$K_2(s) \in C^1([0,\tau_2])$
such that
$$
K_1(s)=K_1^{*}(s)>0, \quad \frac{d}{ds}K_1(s)<0, \quad s \in [0,\tau_1],
\eqno (19)
$$
$$
K_2(s)=K_2^{*}(s)>0, \quad \frac{d}{ds}K_2(s)<0, \quad s \in [0,\tau_2],
\eqno (20)
$$
$$
C=-\begin{pmatrix}
HA+A^{*}H+K_1(0)+K_2(0) && HB_1 && HB_2 \\
B_1^{*}H && -K_1(\tau_1) && 0 \\
B_2^{*}H && 0 && -K_2(\tau_2)
\end{pmatrix}>0,
\eqno (21)
$$
then the zero solution to system~(18) is asymptotically stable.
This result is proved by analogy with~\cite{DM2005}
using the Lyapunov--Krasovskii functional
$$
V(t,{\bf y})=\left\langle H{\bf y}(t),{\bf y}(t) \right\rangle
+\int\limits_{t-\tau_1}^{t}
\left\langle K_1(t-s){\bf y}(s),{\bf y}(s) \right\rangle ds
$$
$$
+\int\limits_{t-\tau_2}^{t}
\left\langle K_2(t-s){\bf y}(s),{\bf y}(s) \right\rangle ds,
\eqno (22)
$$
the derivative of which by virtue of system~(18) has the form
$$
\frac{d}{dt}V(t,{\bf y})
=-\left\langle
C
\begin{pmatrix} {\bf y}(t) \\ {\bf y}(t-\tau_1) \\ {\bf y}(t-\tau_2) \end{pmatrix},
\begin{pmatrix} {\bf y}(t) \\ {\bf y}(t-\tau_1) \\ {\bf y}(t-\tau_2) \end{pmatrix}
\right\rangle
$$
$$
+\int\limits_{t-\tau_1}^{t}
\left\langle \frac{d}{dt}K_1(t-s){\bf y}(s),{\bf y}(s) \right\rangle ds
+\int\limits_{t-\tau_2}^{t}
\left\langle \frac{d}{dt}K_2(t-s){\bf y}(s),{\bf y}(s) \right\rangle ds.
$$

We apply this result in the case when matrices
$A$,
$B_1$,
$B_2$
have the form~(6)--(7).
Our aim is to construct matrices
$H=H^{*}>0$,
$K_1(s) \in C^1([0,\tau_1])$,
$K_2(s) \in C^1([0,\tau_2])$
so that conditions~(19)--(21) are met.

As it was noted above, we assume that condition~(12) is satisfied,
which guarantees the asymptotic stability of equilibrium point
$(x_0,y_0,0)$
of system~(1).
By virtue of formula~(3), condition~(12) is equivalent to the following inequalities:
$$
0<c_1 y_0<\frac{2r}{K} x_0,
\quad
0<e_2 c_2 y_0<d_2.
\eqno (23)
$$

We put
$$
H=\begin{pmatrix}
h_{11} & h_{12} & 0 \\
h_{12} & h_{22} & 0 \\
0 & 0 & h_{33}
\end{pmatrix},
\quad h_{11},h_{22},h_{33}>0, \quad h_{11}h_{22}-h_{12}^2>0,
\eqno (24)
$$
$$
K_1(s)=e^{-m_1 s} (\alpha B_1^{*} B_1+R_1),
\quad m_1>0, \quad \alpha>0,
\quad s \in [0,\tau_1],
\eqno (25)
$$
$$
K_2(s)=e^{-m_2 s} (\beta B_2^{*} B_2+R_2),
\quad m_2>0, \quad \beta>0,
\quad s \in [0,\tau_2],
\eqno (26)
$$
where values
$h_{ij}$,
$m_1$,
$m_2$,
$\alpha$,
$\beta$
and matrices
$R_1=R_1^{*} > 0$,
$R_2=R_2^{*} > 0$
will be defined below.
Taking into account the explicit form~(7) of matrices
$B_1$,
$B_2$
and the explicit form~(24) of matrix
$H$,
we have the following relations:
$B_1^{*}H=B_1^{*}\widetilde{H}_1$,
$B_2^{*}H=B_2^{*}\widetilde{H}_2$,
where
$$
\widetilde{H}_1=\begin{pmatrix}
0 & 0 & 0 \\
h_{12} & h_{22} & 0 \\
0 & 0 & 0
\end{pmatrix},
\quad
\widetilde{H}_2=\begin{pmatrix}
0 & 0 & 0 \\
0 & 0 & 0 \\
0 & 0 & h_{33}
\end{pmatrix}.
$$
Then matrix
$C$
from~(21) will have the form
$$
C=\begin{pmatrix}
-\left(
HA+A^{*}H+\alpha B_1^{*} B_1+\beta B_2^{*} B_2
+\frac{1}{\alpha} e^{m_1 \tau_1} \widetilde{H}_1^{*}\widetilde{H}_1
+\frac{1}{\beta} e^{m_2 \tau_2} \widetilde{H}_2^{*}\widetilde{H}_2
\right) && 0 && 0 \\
0 && 0 && 0 \\
0 && 0 && 0
\end{pmatrix}
$$
$$
+\begin{pmatrix}
\frac{1}{\alpha} e^{m_1 \tau_1} \widetilde{H}_1^{*}\widetilde{H}_1 &&
-\widetilde{H}_1^{*}B_1 && 0 \\
-B_1^{*}\widetilde{H}_1 && \alpha e^{-m_1 \tau_1} B_1^{*} B_1 && 0 \\
0 && 0 && 0
\end{pmatrix}
+\begin{pmatrix}
\frac{1}{\beta} e^{m_2 \tau_2} \widetilde{H}_2^{*}\widetilde{H}_2 &&
0 && -\widetilde{H}_2^{*}B_2 \\
0 && 0 && 0 \\
-B_2^{*}\widetilde{H}_2 && 0 && \beta e^{-m_2 \tau_2} B_2^{*} B_2
\end{pmatrix}
$$
$$
+\begin{pmatrix}
-R_1-R_2 && 0 && 0 \\
0 && e^{-m_1 \tau_1} R_1 && 0 \\
0 && 0 && e^{-m_2 \tau_2} R_2
\end{pmatrix},
$$
from here
$$
C \ge \begin{pmatrix}
L-R_1-R_2 && 0 && 0 \\
0 && e^{-m_1 \tau_1} R_1 && 0 \\
0 && 0 && e^{-m_2 \tau_2} R_2
\end{pmatrix},
\eqno (27)
$$
where
$$
L=-\left(
HA+A^{*}H+\alpha B_1^{*} B_1+\beta B_2^{*} B_2
+\frac{1}{\alpha} e^{m_1 \tau_1} \widetilde{H}_1^{*}\widetilde{H}_1
+\frac{1}{\beta} e^{m_2 \tau_2} \widetilde{H}_2^{*}\widetilde{H}_2
\right).
\eqno (28)
$$
We have reduced our task to checking the fulfillment of the condition
$L>0$.
For
$\|R_1\| \ll 1$,
$\|R_2\| \ll 1$,
this will lead to the positive definiteness of matrix
$C$.

Taking into account the explicit form of matrices
$A$,
$B_1$,
$B_2$,
$H$,
$\widetilde{H}_1$,
$\widetilde{H}_2$,
we calculate the elements of matrix
$L=(l_{ij})$.
We have
$$
\left\{
\begin{array}{l}
\displaystyle
l_{11}=2h_{11} \frac{r}{K} x_0
-\alpha (e_1 c_1 y_0)^2
-\frac{1}{\alpha} e^{m_1 \tau_1} h_{12}^2,
\\
\\
\displaystyle
l_{12}=h_{11} \frac{d_1}{e_1}
+h_{12} \left( \frac{r}{K} x_0+d_1 \right)
-\alpha (e_1 c_1 y_0) d_1
-\frac{1}{\alpha} e^{m_1 \tau_1} h_{12}h_{22},
\\
\\
\displaystyle
l_{22}=2h_{12} \frac{d_1}{e_1}+2h_{22} d_1
-\alpha d_1^2
-\frac{1}{\alpha} e^{m_1 \tau_1} h_{22}^2,
\\
\\
\displaystyle
l_{13}=h_{12} c_2 y_0,
\\
\\
\displaystyle
l_{23}=h_{22} c_2 y_0,
\\
\\
\displaystyle
l_{33}=2h_{33} d_2
-\beta (e_2 c_2 y_0)^2
-\frac{1}{\beta} e^{m_2 \tau_2} h_{33}^2.
\end{array}
\right.
$$
Assuming
$$
\beta=h_{33} \frac{e^{m_2 \tau_2 /2}}{e_2 c_2 y_0},
\eqno (29)
$$
we obtain
$$
l_{33}=2h_{33} \left( d_2-e_2 c_2 y_0 e^{m_2 \tau_2 /2} \right).
$$
Inequality
$l_{33}>0$
will be executed if the number
$m_2>0$
satisfies the condition
$$
e_2 c_2 y_0 e^{m_2 \tau_2 /2} < d_2.
\eqno (30)
$$
Due to the second inequality in formula~(23), such number
$m_2>0$
exists.

Now we check the fulfillment of the conditions
$l_{11}>0$,
$l_{22}>0$,
$l_{11}l_{22}-l_{12}^2>0$.
Assuming
$$
h_{22}=\alpha \left( \frac{r}{K} x_0+d_1 \right) e^{-m_1 \tau_1},
\eqno (31)
$$
we have
$$
\left\{
\begin{array}{l}
\displaystyle
l_{11}=2h_{11} \frac{r}{K} x_0
-\alpha (e_1 c_1 y_0)^2
-\frac{1}{\alpha} e^{m_1 \tau_1} h_{12}^2,
\\
\\
\displaystyle
l_{12}=h_{11} \frac{d_1}{e_1}
-\alpha (e_1 c_1 y_0) d_1,
\\
\\
\displaystyle
l_{22}=2h_{12} \frac{d_1}{e_1}
-\alpha d_1^2 (1-e^{-m_1 \tau_1})
-\alpha \left( \frac{r}{K} x_0 \right)^2 e^{-m_1 \tau_1}.
\end{array}
\right.
$$
We write down the value
$l_{11}l_{22}-l_{12}^2$:
$$
l_{11}l_{22}-l_{12}^2
=\left( 2h_{11} \frac{r}{K} x_0
-\alpha (e_1 c_1 y_0)^2
-\frac{1}{\alpha} e^{m_1 \tau_1} h_{12}^2 \right)
l_{22}
-\left( h_{11} \frac{d_1}{e_1}
-\alpha (e_1 c_1 y_0) d_1 \right)^2
$$
$$
=2h_{11} \left(
\frac{r}{K} x_0 l_{22}+\alpha (c_1 y_0) d_1^2
\right)
-h_{11}^2 \left( \frac{d_1}{e_1} \right)^2
$$
$$
-\alpha (e_1 c_1 y_0)^2 l_{22}
-\frac{1}{\alpha} e^{m_1 \tau_1} h_{12}^2 l_{22}
-\alpha^2 (e_1 c_1 y_0)^2 d_1^2.
$$
We put
$$
h_{11}=\left( \frac{e_1}{d_1} \right)^2
\left(
\frac{r}{K} x_0 l_{22}+\alpha (c_1 y_0) d_1^2
\right),
\eqno (32)
$$
then
$$
l_{11}l_{22}-l_{12}^2
=\left[
\left( \frac{e_1}{d_1} \frac{r}{K} x_0 \right)^2 l_{22}
+2\alpha \frac{r}{K} x_0 (e_1 c_1 y_0) e_1
-\alpha (e_1 c_1 y_0)^2
-\frac{1}{\alpha} e^{m_1 \tau_1} h_{12}^2
\right]
l_{22}
$$
$$
=\bigg[
\alpha e_1^2
\left\{
\left( \frac{r}{K} x_0 \right)^2 e^{-m_1 \tau_1}
-\left( \frac{r}{K} x_0-c_1 y_0 \right)^2
\right\}
$$
$$
-\frac{1}{\alpha} e^{m_1 \tau_1}
\left(
h_{12}
-\alpha \frac{e_1}{d_1} \left( \frac{r}{K} x_0 \right)^2 e^{-m_1 \tau_1}
\right)^2
\bigg]
l_{22}.
$$
Assuming
$$
h_{12}=\alpha \frac{e_1}{d_1}
\left( \frac{r}{K} x_0 \right)^2 e^{-m_1 \tau_1},
\eqno (33)
$$
we obtain
$$
l_{11}l_{22}-l_{12}^2
=\alpha e_1^2
\left\{
\left( \frac{r}{K} x_0 \right)^2 e^{-m_1 \tau_1}
-\left( \frac{r}{K} x_0-c_1 y_0 \right)^2
\right\}
l_{22},
$$
where
$$
\left\{
\begin{array}{l}
\displaystyle
l_{11}=\left( \frac{e_1}{d_1} \frac{r}{K} x_0 \right)^2 l_{22}
+\alpha e_1^2
\left\{
\left( \frac{r}{K} x_0 \right)^2 e^{-m_1 \tau_1}
-\left( \frac{r}{K} x_0-c_1 y_0 \right)^2
\right\},
\\
\\
\displaystyle
l_{12}=\frac{e_1}{d_1} \frac{r}{K} x_0 l_{22},
\\
\\
\displaystyle
l_{22}=\alpha
\left\{
\left( \left( \frac{r}{K} x_0 \right)^2+d_1^2 \right) e^{-m_1 \tau_1}-d_1^2
\right\}.
\end{array}
\right.
$$

We define the number
$m_1>0$
from the conditions:
$$
\left\{
\begin{array}{l}
\displaystyle
\left( \frac{r}{K} x_0 \right)^2 e^{-m_1 \tau_1}
>\left( \frac{r}{K} x_0-c_1 y_0 \right)^2,
\\
\\
\displaystyle
\left( \left( \frac{r}{K} x_0 \right)^2+d_1^2 \right) e^{-m_1 \tau_1}
>d_1^2.
\end{array}
\right.
\eqno (34)
$$
Since the first inequality in formula~(23) is fulfilled,
then such number
$m_1>0$
exists.
Under these conditions on
$m_1$,
the inequalities are valid:
$l_{11}>0$,
$l_{22}>0$,
$l_{11}l_{22}-l_{12}^2>0$.

So, taking into account the above formulas, the elements of matrix
$L$
will have the form:
$$
\left\{
\begin{array}{l}
\displaystyle
l_{11}=\left( \frac{e_1}{d_1} \frac{r}{K} x_0 \right)^2 l_{22}
+\alpha e_1^2
\left\{
\left( \frac{r}{K} x_0 \right)^2 e^{-m_1 \tau_1}
-\left( \frac{r}{K} x_0-c_1 y_0 \right)^2
\right\},
\\
\\
\displaystyle
l_{12}=\frac{e_1}{d_1} \frac{r}{K} x_0 l_{22},
\\
\\
\displaystyle
l_{22}=\alpha
\left\{
\left( \left( \frac{r}{K} x_0 \right)^2+d_1^2 \right) e^{-m_1 \tau_1}-d_1^2
\right\},
\\
\\
\displaystyle
l_{13}=\alpha c_2 y_0 \frac{e_1}{d_1}
\left( \frac{r}{K} x_0 \right)^2 e^{-m_1 \tau_1},
\\
\\
\displaystyle
l_{23}=\alpha c_2 y_0 \left( \frac{r}{K} x_0+d_1 \right) e^{-m_1 \tau_1},
\\
\\
\displaystyle
l_{33}=2h_{33} \left( d_2-e_2 c_2 y_0 e^{m_2 \tau_2 /2} \right).
\end{array}
\right.
\eqno (35)
$$
We choose the number
$h_{33}>0$
from the condition
$\det L>0$.
We have
$$
\det L=(l_{11}l_{22}-l_{12}^2)
\left(
l_{33}-\frac{(l_{11}l_{23}^2+l_{22}l_{13}^2-2l_{12}l_{13}l_{23})}
{(l_{11}l_{22}-l_{12}^2)}
\right)>0,
$$
from here
$$
2h_{33} \left( d_2-e_2 c_2 y_0 e^{m_2 \tau_2 /2} \right)
>\frac{(l_{11}l_{23}^2+l_{22}l_{13}^2-2l_{12}l_{13}l_{23})}
{(l_{11}l_{22}-l_{12}^2)}.
\eqno (36)
$$
Thus, the inequalities are fair:
$l_{11}>0$,
$l_{11}l_{22}-l_{12}^2>0$,
$\det L>0$.
According to Sylvester's criterion, this means that matrix
$L$
is positive definite.

{\bf Remark.}
The number
$\alpha>0$
can be specified arbitrarily, for example, we can take
$\alpha=1$.

It remains to check that
$H>0$.
Indeed, it follows from formula~(28) that
$$
HA+A^{*}H=-\left(
L+\alpha B_1^{*} B_1+\beta B_2^{*} B_2
+\frac{1}{\alpha} e^{m_1 \tau_1} \widetilde{H}_1^{*}\widetilde{H}_1
+\frac{1}{\beta} e^{m_2 \tau_2} \widetilde{H}_2^{*}\widetilde{H}_2
\right)<0,
$$
i.e., matrix
$H=H^{*}$
is a solution to matrix Lyapunov equation
$HA+A^{*}H=-S$,
where
$S=S^{*}>0$.
Since all eigenvalues of matrix
$A$
are contained in the left half-plane
$\mathbb{C}_{-}$,
hence matrix
$H$
is positive definite
(see, for example,~\cite{D2009_2}, chapter~1, section~4).

So, the Lyapunov--Krasovskii functional has been constructed.


\begin{center}
{\bf 4. Obtaining estimates for solutions}
\end{center}

In this section, we indicate conditions for the initial data~(2) for system~(1),
under which solutions stabilize to equilibrium point
$(x_0,y_0,0)$,
and we obtain estimates for solutions characterizing the stabilization rate at infinity.
As above, assume that inequa\-lities~(23) are fulfilled,
which guarantee the asymptotic stability of this equilibrium point.

Change of variables~(4) leads system~(1) to the form~(5),
initial conditions~(2) are rewritten as follows
$$
\left\{
\begin{array}{l}
\displaystyle
\widetilde{x}(\theta)=\varphi(\theta)-x_0 \ge -x_0, \quad \theta \in [-\tau_1,0], \quad \widetilde{x}(+0)=\varphi(0)-x_0>-x_0,
\\
\\
\displaystyle
\widetilde{y}(\theta)=\psi(\theta)-y_0 \ge -y_0, \quad \theta \in [-\tau_{\max},0], \quad \widetilde{y}(+0)=\psi(0)-y_0,
\\
\\
\displaystyle
\widetilde{z}(\theta)=\eta(\theta) \ge 0, \quad \theta \in [-\tau_2,0], \quad \widetilde{z}(+0)=\eta(0).
\end{array}
\right.
\eqno (37)
$$

Below we will formulate the main result of this section.
To do this, we first introduce the notation.

Define vector function
$\widetilde{\bm\psi}(\theta)$
by the following rule:
$$
\widetilde{\bm \psi}(\theta)
=\begin{pmatrix}
\widetilde{\varphi}(\theta) \\
\widetilde{\psi}(\theta) \\
\widetilde{\eta}(\theta)
\end{pmatrix},
\quad t \leq 0,
$$
where
$$
\widetilde{\varphi}(\theta)
=\left\{
\begin{array}{l}
\varphi(\theta)-x_0,
\\
0,
\end{array}
\quad
\begin{array}{l}
\theta \in [-\tau_1,0],
\\
\theta < -\tau_1,
\end{array}
\right.
$$
$$
\widetilde{\psi}(\theta)
=\left\{
\begin{array}{l}
\psi(\theta)-y_0,
\\
0,
\end{array}
\quad
\begin{array}{l}
\theta \in [-\tau_{\max},0],
\\
\theta < -\tau_{\max},
\end{array}
\right.
$$
$$
\widetilde{\eta}(\theta)
=\left\{
\begin{array}{l}
\eta(\theta),
\\
0,
\end{array}
\quad
\begin{array}{l}
\theta \in [-\tau_2,0],
\\
\theta < -\tau_2.
\end{array}
\right.
$$
Consider Lyapunov--Krasovskii functional of the form~(22)
$$
V(0,\widetilde{\bm \psi})
=\left\langle H\widetilde{\bm \psi}(0),\widetilde{\bm \psi}(0) \right\rangle
+\int\limits_{-\tau_1}^{0}
\left\langle K_1(-\theta)\widetilde{\bm \psi}(\theta),
\widetilde{\bm \psi}(\theta) \right\rangle d\theta
+\int\limits_{-\tau_2}^{0}
\left\langle K_2(-\theta)\widetilde{\bm \psi}(\theta),
\widetilde{\bm \psi}(\theta) \right\rangle d\theta,
$$
where matrix
$H=H^{*}>0$
has the form~(24), its elements
$h_{ij}$
are defined in (31), (32), (33), (36); matrices
$K_1(s)$,
$K_2(s)$
have the form~(25), (26):
$$
K_1(s)=e^{-m_1 s} (\alpha B_1^{*} B_1+R_1),
\quad s \in [0,\tau_1],
\eqno (38)
$$
$$
K_2(s)=e^{-m_2 s} (\beta B_2^{*} B_2+R_2),
\quad s \in [0,\tau_2],
\eqno (39)
$$
$\alpha=1$,
$\beta>0$
is defined in~(29),
$m_1>0$
is defined in~(34),
$m_2>0$
is defined in~(30),
$R_1=R_1^{*} \ge 0$
and
$R_2=R_2^{*} \ge 0$
will be defined below.

Let
$\sigma>0$
be such a number that the inequality holds
$$
L \ge \sigma H,
\eqno (40)
$$
where matrix
$L=L^{*}>0$
is defined in~(28), its elements have the form~(35).
For example, we can take
$\sigma=\lambda_{\min}(H^{-1/2} L H^{-1/2})$,
the minimal eigenvalue of matrix
$H^{-1/2} L H^{-1/2}$.

Next, let values
$\mu_1>0$
and
$\mu_2>0$
satisfy the inequality
$$
\max\{ \mu_1, \mu_2 \}<\frac{\sigma}{2}.
$$
We put
$$
R_1=\mu_1 H_1,
\quad
H_1=\begin{pmatrix}
h_{11} & h_{12} & 0 \\
h_{12} & h_{22} & 0 \\
0 & 0 & 0
\end{pmatrix},
\eqno (41)
$$
$$
R_2=\mu_2 H_2,
\quad
H_2=\begin{pmatrix}
0 & 0 & 0 \\
0 & 0 & 0 \\
0 & 0 & h_{33}
\end{pmatrix}.
\eqno (42)
$$

Finally, we denote
$$
\varepsilon=\min\left\{
\Big( \sigma-2\max\{ \mu_1, \mu_2 \} \Big), m_1, m_2
\right\},
\eqno (43)
$$
$$
q=\frac{2}{\sqrt{\left( 1-\frac{h_{12}}{\sqrt{h_{11}h_{22}}} \right)}}
\max\left\{
\frac{\sqrt{\left( \frac{r}{K} \right)^2+c_1^2}}
{\min\left\{ \sqrt{h_{11}}, \sqrt{h_{22}} \right\}
\sqrt{\left( 1-\frac{h_{12}}{\sqrt{h_{11}h_{22}}} \right)}}, \
\frac{c_2}{\sqrt{h_{33}}}
\right\}.
\eqno (44)
$$

The following theorem is valid.

{\bf Theorem.}
{\it
Let inequalities~$(23)$ be fulfilled and initial functions
$\varphi(\theta)$,
$\psi(\theta)$,
$\eta(\theta)$
satisfy the conditions:
$$
\max\limits_{\theta \in [-\tau_1,0]} |\psi(\theta)-y_0|
\leq \frac{\sqrt{h_{11}h_{22}-h_{12}^2}}{h_{22}}
\left( \frac{\mu_1}{e_1 c_1} \right) e^{-m_1 \tau_1 /2},
\eqno (45)
$$
$$
\max\limits_{\theta \in [-\tau_2,0]} |\psi(\theta)-y_0|
\leq \left( \frac{\mu_2}{e_2 c_2} \right) e^{-m_2 \tau_2 /2},
\eqno (46)
$$
$$
\sqrt{V(0,\widetilde{\bm \psi})}<\frac{\varepsilon}{q},
\eqno (47)
$$
$$
\frac{\sqrt{V(0,\widetilde{\bm \psi})}}
{\displaystyle
\left( 1-\frac{q}{\varepsilon} \sqrt{V(0,\widetilde{\bm \psi})} \right)}
\leq \frac{\left( h_{11}h_{22}-h_{12}^2 \right)}{h_{22} \sqrt{h_{11}}}
\left( \frac{\mu_1}{e_1 c_1} \right) e^{-m_1 \tau_1 /2},
\eqno (48)
$$
$$
\frac{\sqrt{V(0,\widetilde{\bm \psi})}}
{\displaystyle
\left( 1-\frac{q}{\varepsilon} \sqrt{V(0,\widetilde{\bm \psi})} \right)}
\leq \frac{\sqrt{h_{11}h_{22}-h_{12}^2}}{\sqrt{h_{11}}}
\left( \frac{\mu_2}{e_2 c_2} \right) e^{-m_2 \tau_2 /2}.
\eqno (49)
$$
Then for the components of the solution to the initial value problem~$(1)$, $(2)$,
the estimates are fair
$$
|x(t)-x_0| \leq \frac{\sqrt{h_{22}}}{\sqrt{h_{11}h_{22}-h_{12}^2}}
\frac{\sqrt{V(0,\widetilde{\bm \psi})} e^{-\varepsilon t/2}}
{\displaystyle
\left( 1-\frac{q}{\varepsilon} \sqrt{V(0,\widetilde{\bm \psi})} \right)},
\quad t>0,
\eqno (50)
$$
$$
|y(t)-y_0| \leq \frac{\sqrt{h_{11}}}{\sqrt{h_{11}h_{22}-h_{12}^2}}
\frac{\sqrt{V(0,\widetilde{\bm \psi})} e^{-\varepsilon t/2}}
{\displaystyle
\left( 1-\frac{q}{\varepsilon} \sqrt{V(0,\widetilde{\bm \psi})} \right)},
\quad t>0,
\eqno (51)
$$
$$
|z(t)| \leq \frac{1}{\sqrt{h_{33}}}
\frac{\sqrt{V(0,\widetilde{\bm \psi})} e^{-\varepsilon t/2}}
{\displaystyle
\left( 1-\frac{q}{\varepsilon} \sqrt{V(0,\widetilde{\bm \psi})} \right)},
\quad t>0.
\eqno (52)
$$
}

{\bf Proof.}
1. Let
$\widetilde{\bf y}(t)$
be a solution to the initial value problem~(5), (37).
Consider the Lyapunov--Krasovskii functional of the form~(22):
$$
V(t, \widetilde{\bf y})
=\left\langle H\widetilde{\bf y}(t), \widetilde{\bf y}(t) \right\rangle
+\int\limits_{t-\tau_1}^{t}
\left\langle K_1(t-s) \widetilde{\bf y}(s), \widetilde{\bf y}(s) \right\rangle
ds
$$
$$
+\int\limits_{t-\tau_2}^{t}
\left\langle K_2(t-s) \widetilde{\bf y}(s), \widetilde{\bf y}(s) \right\rangle
ds.
\eqno (53)
$$
Its derivative by virtue of system~(5) has the form
$$
\frac{d}{dt}V(t,\widetilde{\bf y})
=-\left\langle
C
\begin{pmatrix}
\widetilde{\bf y}(t) \\
\widetilde{\bf y}(t-\tau_1) \\
\widetilde{\bf y}(t-\tau_2)
\end{pmatrix},
\begin{pmatrix}
\widetilde{\bf y}(t) \\
\widetilde{\bf y}(t-\tau_1) \\
\widetilde{\bf y}(t-\tau_2)
\end{pmatrix}
\right\rangle
$$
$$
+2\left\langle H\widetilde{\bf y}(t), F(\widetilde{\bf y}(t)) \right\rangle
+2\left\langle H\widetilde{\bf y}(t), G_1(\widetilde{\bf y}(t-\tau_1)) \right\rangle
+2\left\langle H\widetilde{\bf y}(t), G_2(\widetilde{\bf y}(t-\tau_2)) \right\rangle
$$
$$
+\int\limits_{t-\tau_1}^{t}
\left\langle \frac{d}{dt}K_1(t-s)\widetilde{\bf y}(s),
\widetilde{\bf y}(s) \right\rangle ds
+\int\limits_{t-\tau_2}^{t}
\left\langle \frac{d}{dt}K_2(t-s)\widetilde{\bf y}(s),
\widetilde{\bf y}(s) \right\rangle ds,
$$
where matrix
$C$
is defined in~(21).
From inequalities~(27) and~(40), the estimate follows
$$
-\left\langle
C
\begin{pmatrix}
\widetilde{\bf y}(t) \\
\widetilde{\bf y}(t-\tau_1) \\
\widetilde{\bf y}(t-\tau_2)
\end{pmatrix},
\begin{pmatrix}
\widetilde{\bf y}(t) \\
\widetilde{\bf y}(t-\tau_1) \\
\widetilde{\bf y}(t-\tau_2)
\end{pmatrix}
\right\rangle
\leq -\left\langle (L-R_1-R_2)\widetilde{\bf y}(t),\widetilde{\bf y}(t) \right\rangle
$$
$$
-e^{-m_1 \tau_1}
\left\langle R_1\widetilde{\bf y}(t-\tau_1),\widetilde{\bf y}(t-\tau_1) \right\rangle
-e^{-m_2 \tau_2}
\left\langle R_2\widetilde{\bf y}(t-\tau_2),\widetilde{\bf y}(t-\tau_2) \right\rangle
$$
$$
\leq -\sigma
\left\langle H\widetilde{\bf y}(t),\widetilde{\bf y}(t) \right\rangle
+\left\langle R_1\widetilde{\bf y}(t),\widetilde{\bf y}(t) \right\rangle
+\left\langle R_2\widetilde{\bf y}(t),\widetilde{\bf y}(t) \right\rangle
$$
$$
-e^{-m_1 \tau_1}
\left\langle R_1\widetilde{\bf y}(t-\tau_1),\widetilde{\bf y}(t-\tau_1) \right\rangle
-e^{-m_2 \tau_2}
\left\langle R_2\widetilde{\bf y}(t-\tau_2),\widetilde{\bf y}(t-\tau_2) \right\rangle,
$$
and from formulas~(38)--(39), the relations follow
$$
\frac{d}{dt}K_1(t-s)=-m_1 K_1(t-s),
\quad s \in [t-\tau_1,t],
$$
$$
\frac{d}{dt}K_2(t-s)=-m_2 K_2(t-s),
\quad s \in [t-\tau_2,t].
$$
Then the inequality is fair
$$
\frac{d}{dt}V(t,\widetilde{\bf y})
\leq -\sigma
\left\langle H\widetilde{\bf y}(t),\widetilde{\bf y}(t) \right\rangle
+g_1(t)
+g_2(t)
+f(t)
$$
$$
-m_1 \int\limits_{t-\tau_1}^{t}
\left\langle K_1(t-s)\widetilde{\bf y}(s),
\widetilde{\bf y}(s) \right\rangle ds
-m_2 \int\limits_{t-\tau_2}^{t}
\left\langle K_2(t-s)\widetilde{\bf y}(s),
\widetilde{\bf y}(s) \right\rangle ds,
\eqno (54)
$$
where
$$
g_1(t)=\left\langle R_1\widetilde{\bf y}(t),\widetilde{\bf y}(t) \right\rangle
-e^{-m_1 \tau_1}
\left\langle
R_1
\widetilde{\bf y}(t-\tau_1),
\widetilde{\bf y}(t-\tau_1)
\right\rangle
+2\left\langle H\widetilde{\bf y}(t), G_1(\widetilde{\bf y}(t-\tau_1)) \right\rangle,
$$
$$
g_2(t)=\left\langle R_2\widetilde{\bf y}(t),\widetilde{\bf y}(t) \right\rangle
-e^{-m_2 \tau_2}
\left\langle
R_2
\widetilde{\bf y}(t-\tau_2),
\widetilde{\bf y}(t-\tau_2)
\right\rangle
+2\left\langle H\widetilde{\bf y}(t), G_2(\widetilde{\bf y}(t-\tau_2)) \right\rangle,
$$
$$
f(t)=2\left\langle H\widetilde{\bf y}(t), F(\widetilde{\bf y}(t)) \right\rangle.
$$

2. We will evaluate
$g_1(t)$.
Taking into account the explicit form of matrix
$R_1$
(see~(41))
and vector function
$G_1(\widetilde{\bf y}(t-\tau_1))$
(see~(9)),
we get
$$
g_1(t)=\mu_1 \left\langle H_1 \widetilde{\bf y}(t), \widetilde{\bf y}(t) \right\rangle
-\mu_1 e^{-m_1 \tau_1} \left\langle H_1 \widetilde{\bf y}(t-\tau_1), \widetilde{\bf y}(t-\tau_1) \right\rangle
+2\left\langle H_1\widetilde{\bf y}(t), G_1(\widetilde{\bf y}(t-\tau_1)) \right\rangle
$$
$$
\leq \mu_1 \left\langle H_1 \widetilde{\bf y}(t), \widetilde{\bf y}(t) \right\rangle
-\mu_1 e^{-m_1 \tau_1} \left\langle H_1 \widetilde{\bf y}(t-\tau_1), \widetilde{\bf y}(t-\tau_1) \right\rangle
$$
$$
+2\sqrt{\left\langle H_1\widetilde{\bf y}(t), \widetilde{\bf y}(t) \right\rangle}
\sqrt{\left\langle H_1G_1(\widetilde{\bf y}(t-\tau_1)), G_1(\widetilde{\bf y}(t-\tau_1)) \right\rangle}.
$$
Due to inequality
$$
\sqrt{\left\langle H_1G_1(\widetilde{\bf y}(t-\tau_1)), G_1(\widetilde{\bf y}(t-\tau_1)) \right\rangle}
=e_1 c_1 \sqrt{h_{22}} \
|\widetilde{x}(t-\tau_1)| |\widetilde{y}(t-\tau_1)|
$$
$$
\leq e_1 c_1 \frac{h_{22}}{\sqrt{h_{11}h_{22}-h_{12}^2}} |\widetilde{y}(t-\tau_1)|
\sqrt{\left\langle H_1\widetilde{\bf y}(t-\tau_1), \widetilde{\bf y}(t-\tau_1) \right\rangle},
$$
we establish the estimate for function
$g_1(t)$:
$$
g_1(t) \leq \mu_1 \left\langle H_1 \widetilde{\bf y}(t), \widetilde{\bf y}(t) \right\rangle
-\mu_1 e^{-m_1 \tau_1} \left\langle H_1 \widetilde{\bf y}(t-\tau_1), \widetilde{\bf y}(t-\tau_1) \right\rangle
$$
$$
+2 e_1 c_1 \frac{h_{22}}{\sqrt{h_{11}h_{22}-h_{12}^2}} |\widetilde{y}(t-\tau_1)|
\sqrt{\left\langle H_1\widetilde{\bf y}(t), \widetilde{\bf y}(t) \right\rangle}
\sqrt{\left\langle H_1\widetilde{\bf y}(t-\tau_1), \widetilde{\bf y}(t-\tau_1) \right\rangle}
$$
$$
\leq \mu_1 \left(
1+e^{m_1 \tau_1} \left( \frac{e_1 c_1}{\mu_1} \right)^2
\frac{h_{22}^2}{(h_{11}h_{22}-h_{12}^2)}
|\widetilde{y}(t-\tau_1)|^2
\right)
\left\langle H_1 \widetilde{\bf y}(t), \widetilde{\bf y}(t) \right\rangle.
\eqno (55)
$$

3. We will evaluate
$g_2(t)$.
Taking into account the explicit form of matrix
$R_2$
(see~(42))
and vector function
$G_2(\widetilde{\bf y}(t-\tau_2))$
(see~(10)),
we obtain the estimate
$$
g_2(t)=\mu_2 \left\langle H_2 \widetilde{\bf y}(t), \widetilde{\bf y}(t) \right\rangle
-\mu_2 e^{-m_2 \tau_2} \left\langle H_2 \widetilde{\bf y}(t-\tau_2), \widetilde{\bf y}(t-\tau_2) \right\rangle
+2\left\langle H_2\widetilde{\bf y}(t), G_2(\widetilde{\bf y}(t-\tau_2)) \right\rangle
$$
$$
=h_{33}
\Big(
\mu_2 \widetilde{z}^2(t)
-\mu_2 e^{-m_2 \tau_2} \widetilde{z}^2(t-\tau_2)
+2e_2 c_2 \widetilde{y}(t-\tau_2)
\widetilde{z}(t-\tau_2) \widetilde{z}(t)
\Big)
$$
$$
\leq \mu_2 \left(
1+e^{m_2 \tau_2} \left( \frac{e_2 c_2}{\mu_2} \right)^2
|\widetilde{y}(t-\tau_2)|^2
\right)
h_{33} \widetilde{z}^2(t)
$$
$$
=\mu_2 \left(
1+e^{m_2 \tau_2} \left( \frac{e_2 c_2}{\mu_2} \right)^2
|\widetilde{y}(t-\tau_2)|^2
\right)
\left\langle H_2 \widetilde{\bf y}(t), \widetilde{\bf y}(t) \right\rangle.
\eqno (56)
$$

4. We will evaluate
$f(t)$.
Taking into account definition~(8) of vector function
$F(\widetilde{\bf y}(t))$,
we rewrite function
$f(t)$
in the following form:
$$
f(t)=2\left\langle H_1\widetilde{\bf y}(t), F(\widetilde{\bf y}(t)) \right\rangle,
$$
where matrix
$H_1$
is defined in~(41).
From here, the estimate follows
$$
f(t) \leq 2\sqrt{\left\langle H_1\widetilde{\bf y}(t), \widetilde{\bf y}(t) \right\rangle}
\sqrt{\left\langle H_1F(\widetilde{\bf y}(t)), F(\widetilde{\bf y}(t)) \right\rangle}
=2\sqrt{\left\langle H_1\widetilde{\bf y}(t), \widetilde{\bf y}(t) \right\rangle}
$$
$$
\times
\sqrt{\left\langle
\begin{pmatrix}
h_{11} & h_{12} \\
h_{12} & h_{22}
\end{pmatrix}
\begin{pmatrix}
\widetilde{x}(t) & 0 \\
0 & \widetilde{y}(t)
\end{pmatrix}
\begin{pmatrix}
\frac{r}{K} \widetilde{x}(t)+c_1 \widetilde{y}(t) \\
c_2 \widetilde{z}(t)
\end{pmatrix},
\begin{pmatrix}
\widetilde{x}(t) & 0 \\
0 & \widetilde{y}(t)
\end{pmatrix}
\begin{pmatrix}
\frac{r}{K} \widetilde{x}(t)+c_1 \widetilde{y}(t) \\
c_2 \widetilde{z}(t)
\end{pmatrix}
\right\rangle}
$$
$$
=2\sqrt{\left\langle H_1\widetilde{\bf y}(t), \widetilde{\bf y}(t) \right\rangle}
$$
$$
\times
\sqrt{\left\langle
\begin{pmatrix}
h_{11} \widetilde{x}^2(t) & h_{12} \widetilde{x}(t)\widetilde{y}(t) \\
h_{12} \widetilde{x}(t)\widetilde{y}(t) & h_{22} \widetilde{y}^2(t)
\end{pmatrix}
\begin{pmatrix}
\frac{r}{K} \widetilde{x}(t)+c_1 \widetilde{y}(t) \\
c_2 \widetilde{z}(t)
\end{pmatrix},
\begin{pmatrix}
\frac{r}{K} \widetilde{x}(t)+c_1 \widetilde{y}(t) \\
c_2 \widetilde{z}(t)
\end{pmatrix}
\right\rangle}
$$
$$
\leq 2\gamma(t)\delta(t)
\sqrt{\left\langle H\widetilde{\bf y}(t), \widetilde{\bf y}(t) \right\rangle},
\eqno (57)
$$
where
$$
\gamma(t)=\sqrt{\left\|
\begin{pmatrix}
h_{11} \widetilde{x}^2(t) & h_{12} \widetilde{x}(t)\widetilde{y}(t) \\
h_{12} \widetilde{x}(t)\widetilde{y}(t) & h_{22} \widetilde{y}^2(t)
\end{pmatrix}
\right\|},
$$
$$
\delta(t)=\left\|
\begin{pmatrix}
\frac{r}{K} \widetilde{x}(t)+c_1 \widetilde{y}(t) \\
c_2 \widetilde{z}(t)
\end{pmatrix}
\right\|.
\eqno (58)
$$

First we get an estimate for
$\gamma(t)$.
Due to inequalities
$$
\begin{pmatrix}
h_{11} \widetilde{x}^2(t) & h_{12} \widetilde{x}(t)\widetilde{y}(t) \\
h_{12} \widetilde{x}(t)\widetilde{y}(t) & h_{22} \widetilde{y}^2(t)
\end{pmatrix}
\leq \left( h_{11} \widetilde{x}^2(t)+h_{22} \widetilde{y}^2(t) \right)
\begin{pmatrix}
1 & 0 \\
0 & 1
\end{pmatrix},
$$
$$
\begin{pmatrix}
h_{11} & 0 \\
0 & h_{22}
\end{pmatrix}
\leq \frac{1}{\left( 1-\frac{h_{12}}{\sqrt{h_{11}h_{22}}} \right)}
\begin{pmatrix}
h_{11} & h_{12} \\
h_{12} & h_{22}
\end{pmatrix},
$$
the estimate is fair
$$
\gamma(t) \leq \sqrt{h_{11} \widetilde{x}^2(t)+h_{22} \widetilde{y}^2(t)}
=\sqrt{\left\langle
\begin{pmatrix}
h_{11} & 0 \\
0 & h_{22}
\end{pmatrix}
\begin{pmatrix}
\widetilde{x}(t) \\ \widetilde{y}(t)
\end{pmatrix},
\begin{pmatrix}
\widetilde{x}(t) \\ \widetilde{y}(t)
\end{pmatrix}
\right\rangle}
$$
$$
\leq \frac{1}{\sqrt{\left( 1-\frac{h_{12}}{\sqrt{h_{11}h_{22}}} \right)}}
\sqrt{\left\langle
\begin{pmatrix}
h_{11} & h_{12} \\
h_{12} & h_{22}
\end{pmatrix}
\begin{pmatrix}
\widetilde{x}(t) \\ \widetilde{y}(t)
\end{pmatrix},
\begin{pmatrix}
\widetilde{x}(t) \\ \widetilde{y}(t)
\end{pmatrix}
\right\rangle}
$$
$$
\leq \frac{1}{\sqrt{\left( 1-\frac{h_{12}}{\sqrt{h_{11}h_{22}}} \right)}}
\sqrt{\left\langle H\widetilde{\bf y}(t), \widetilde{\bf y}(t) \right\rangle}.
\eqno (59)
$$

Now we evaluate
$\delta(t)$.
Taking into account definition~(58), function
$\delta(t)$
can be rewritten in the form
$$
\delta(t)=\sqrt{
\left\langle
\begin{pmatrix}
\left( \frac{r}{K} \right)^2 & \frac{r}{K} c_1 & 0 \\
\frac{r}{K} c_1 & c_1^2 & 0 \\
0 & 0 & c_2^2
\end{pmatrix}
\begin{pmatrix}
\widetilde{x}(t) \\ \widetilde{y}(t) \\ \widetilde{z}(t)
\end{pmatrix},
\begin{pmatrix}
\widetilde{x}(t) \\ \widetilde{y}(t) \\ \widetilde{z}(t)
\end{pmatrix}
\right\rangle
}.
$$
Due to inequality
$$
\begin{pmatrix}
\left( \frac{r}{K} \right)^2 & \frac{r}{K} c_1 \\
\frac{r}{K} c_1 & c_1^2
\end{pmatrix}
\leq \Big( \left( \frac{r}{K} \right)^2+c_1^2 \Big)
\begin{pmatrix}
1 & 0 \\
0 & 1
\end{pmatrix}
\leq \frac{\Big( \left( \frac{r}{K} \right)^2+c_1^2 \Big)}
{\min\{ h_{11}, h_{22} \}}
\begin{pmatrix}
h_{11} & 0 \\
0 & h_{22}
\end{pmatrix}
$$
$$
\leq \frac{\Big( \left( \frac{r}{K} \right)^2+c_1^2 \Big)}
{\min\{ h_{11}, h_{22} \}
\left( 1-\frac{h_{12}}{\sqrt{h_{11}h_{22}}} \right)}
\begin{pmatrix}
h_{11} & h_{12} \\
h_{12} & h_{22}
\end{pmatrix},
$$
for function
$\delta(t)$,
the estimate holds
$$
\delta(t)
\leq \sqrt{
\frac{\Big( \left( \frac{r}{K} \right)^2+c_1^2 \Big)}
{\min\{ h_{11}, h_{22} \}
\left( 1-\frac{h_{12}}{\sqrt{h_{11}h_{22}}} \right)}
\left\langle H_1\widetilde{\bf y}(t), \widetilde{\bf y}(t) \right\rangle
+\frac{c_2^2}{h_{33}}
\left\langle H_2\widetilde{\bf y}(t), \widetilde{\bf y}(t) \right\rangle
}
$$
$$
\leq \max\left\{
\frac{\sqrt{\left( \frac{r}{K} \right)^2+c_1^2}}
{\min\left\{ \sqrt{h_{11}}, \sqrt{h_{22}} \right\}
\sqrt{\left( 1-\frac{h_{12}}{\sqrt{h_{11}h_{22}}} \right)}}, \
\frac{c_2}{\sqrt{h_{33}}}
\right\}
\sqrt{\left\langle H\widetilde{\bf y}(t), \widetilde{\bf y}(t) \right\rangle}.
\eqno (60)
$$

So, from inequalities~(57), (59), (60), the estimate follows for function
$f(t)$:
$$
f(t) \leq q \left\langle H\widetilde{\bf y}(t), \widetilde{\bf y}(t) \right\rangle^{3/2}
\leq q V^{3/2}(t, \widetilde{\bf y}),
\eqno (61)
$$
where
$q$
is defined in~(44).

5. By virtue of estimates~(55), (56), (61), from inequality~(54), we get
$$
\frac{d}{dt}V(t,\widetilde{\bf y})
\leq -\sigma \left\langle H\widetilde{\bf y}(t),\widetilde{\bf y}(t) \right\rangle
+2\mu_1 \left\langle H_1 \widetilde{\bf y}(t), \widetilde{\bf y}(t) \right\rangle
+2\mu_2 \left\langle H_2 \widetilde{\bf y}(t), \widetilde{\bf y}(t) \right\rangle
$$
$$
-\mu_1 \left(
1-e^{m_1 \tau_1} \left( \frac{e_1 c_1}{\mu_1} \right)^2
\frac{h_{22}^2}{(h_{11}h_{22}-h_{12}^2)}
|\widetilde{y}(t-\tau_1)|^2
\right)
\left\langle H_1 \widetilde{\bf y}(t), \widetilde{\bf y}(t) \right\rangle
$$
$$
-\mu_2 \left(
1-e^{m_2 \tau_2} \left( \frac{e_2 c_2}{\mu_2} \right)^2
|\widetilde{y}(t-\tau_2)|^2
\right)
\left\langle H_2 \widetilde{\bf y}(t), \widetilde{\bf y}(t) \right\rangle
+q V^{3/2}(t, \widetilde{\bf y})
$$
$$
-m_1 \int\limits_{t-\tau_1}^{t}
\left\langle K_1(t-s)\widetilde{\bf y}(s),
\widetilde{\bf y}(s) \right\rangle ds
-m_2 \int\limits_{t-\tau_2}^{t}
\left\langle K_2(t-s)\widetilde{\bf y}(s),
\widetilde{\bf y}(s) \right\rangle ds.
$$

Since
$H=H_1+H_2$,
then
$$
-\sigma \left\langle H\widetilde{\bf y}(t),\widetilde{\bf y}(t) \right\rangle
+2\mu_1 \left\langle H_1 \widetilde{\bf y}(t), \widetilde{\bf y}(t) \right\rangle
+2\mu_2 \left\langle H_2 \widetilde{\bf y}(t), \widetilde{\bf y}(t) \right\rangle
$$
$$
\leq -\Big( \sigma-2\max\{ \mu_1, \mu_2 \} \Big)
\left\langle H \widetilde{\bf y}(t), \widetilde{\bf y}(t) \right\rangle.
$$
Taking into account notation~(43) of value
$\varepsilon$
and definition~(53) of functional
$V(t,\widetilde{\bf y})$,
we establish the inequality
$$
\frac{d}{dt}V(t,\widetilde{\bf y})
\leq -\varepsilon V(t, \widetilde{\bf y})
+q V^{3/2}(t, \widetilde{\bf y})
$$
$$
-\mu_1 \left(
1-e^{m_1 \tau_1} \left( \frac{e_1 c_1}{\mu_1} \right)^2
\frac{h_{22}^2}{(h_{11}h_{22}-h_{12}^2)}
|\widetilde{y}(t-\tau_1)|^2
\right)
\left\langle H_1 \widetilde{\bf y}(t), \widetilde{\bf y}(t) \right\rangle
$$
$$
-\mu_2 \left(
1-e^{m_2 \tau_2} \left( \frac{e_2 c_2}{\mu_2} \right)^2
|\widetilde{y}(t-\tau_2)|^2
\right)
\left\langle H_2 \widetilde{\bf y}(t), \widetilde{\bf y}(t) \right\rangle.
\eqno (62)
$$

6. First, we assume that
$t \in (0,\tau_{\min}]$,
where
$\tau_{\min}=\min\{ \tau_1,\tau_2 \}>0$.
Then
$$
t-\tau_k \leq 0,
\quad k=1,2,
$$
so
$$
\widetilde{y}(t-\tau_k)=\widetilde{\psi}(t-\tau_k)
=\psi(t-\tau_k)-y_0,
\quad k=1,2.
$$
By virtue of conditions~(45) and~(46), from inequality~(62), we have the estimate
$$
\frac{d}{dt}V(t,\widetilde{\bf y})
\leq -\varepsilon V(t, \widetilde{\bf y})
+q V^{3/2}(t, \widetilde{\bf y}).
\eqno (63)
$$
Because by virtue of condition~(47), we have
$\displaystyle
\sqrt{V(0,\widetilde{\bm \psi})}<\frac{\varepsilon}{q}$,
then using the Gronwall inequality
(see, for example,~\cite{Hartman1970}), for
$t \in (0,\tau_{\min}]$,
we get
$$
V(t,\widetilde{\bf y})
\leq \frac{V(0,\widetilde{\bm \psi}) e^{-\varepsilon t}}
{\displaystyle
\left( 1-\frac{q}{\varepsilon} \sqrt{V(0,\widetilde{\bm \psi})} \right)^2}.
$$
Taking into account the estimates
$$
(x(t)-x_0)^2=\widetilde{x}^2(t) \leq \frac{h_{22}}{(h_{11}h_{22}-h_{12}^2)}
\left\langle H \widetilde{\bf y}(t), \widetilde{\bf y}(t) \right\rangle
\leq \frac{h_{22}}{(h_{11}h_{22}-h_{12}^2)} V(t,\widetilde{\bf y}),
$$
$$
(y(t)-y_0)^2=\widetilde{y}^2(t) \leq \frac{h_{11}}{(h_{11}h_{22}-h_{12}^2)}
\left\langle H \widetilde{\bf y}(t), \widetilde{\bf y}(t) \right\rangle
\leq \frac{h_{11}}{(h_{11}h_{22}-h_{12}^2)} V(t,\widetilde{\bf y}),
$$
$$
z^2(t)=\widetilde{z}^2(t) \leq \frac{1}{h_{33}}
\left\langle H \widetilde{\bf y}(t), \widetilde{\bf y}(t) \right\rangle
\leq \frac{1}{h_{33}} V(t,\widetilde{\bf y}),
$$
for
$t \in (0,\tau_{\min}]$,
we establish inequalities~(50)--(52).

Now let
$t \in (\tau_{\min},2\tau_{\min}]$.
Then
$t-\tau_k \leq \tau_{\min}$,
$k=1,2$.
If
$t-\tau_k \leq 0$,
then
$$
\widetilde{y}(t-\tau_k)=\widetilde{\psi}(t-\tau_k)
=\psi(t-\tau_k)-y_0,
\quad k=1,2.
$$
If
$t-\tau_k \in (0,\tau_{\min}]$,
then for function
$\widetilde{y}(t-\tau_k)=(y(t-\tau_k)-y_0)$,
estimate~(51) is fair:
$$
|\widetilde{y}(t-\tau_k)| \leq \frac{\sqrt{h_{11}}}{\sqrt{h_{11}h_{22}-h_{12}^2}}
\frac{\sqrt{V(0,\widetilde{\bm \psi})} e^{-\varepsilon (t-\tau_k)/2}}
{\displaystyle
\left( 1-\frac{q}{\varepsilon} \sqrt{V(0,\widetilde{\bm \psi})} \right)}
$$
$$
\leq \frac{\sqrt{h_{11}}}{\sqrt{h_{11}h_{22}-h_{12}^2}}
\frac{\sqrt{V(0,\widetilde{\bm \psi})}}
{\displaystyle
\left( 1-\frac{q}{\varepsilon} \sqrt{V(0,\widetilde{\bm \psi})} \right)},
\quad k=1,2.
$$
Then by virtue of conditions~(45), (46), (48), (49), from inequality~(62),
we obtain estimate~(63).
Using the same reasoning as in the previous case,
we establish the validity of inequalities~(50)--(52) for
$t \in (\tau_{\min},2\tau_{\min}]$.

Using the method of mathematical induction, for
$t \in (m\tau_{\min},(m+1)\tau_{\min}]$,
$m \in \mathbb{N}$,
the validity of inequalities~(50)--(52)
is proved according to a similar scheme.

So, inequalities~(50)--(52) are proved for all
$t>0$.

Theorem is proved.

{\bf Remark.}
Inequalities~(45)--(49) show at what initial amount
of phytoplankton, zoo\-plankton and initial number of fish,
the plankton survival and fish extinction occur.
Estimates~(50), (51) characterize the stabilization rate
of the amount of phytoplankton and zooplankton to constant values,
estimate~(52) characterizes the decrease rate of the number of fish,
while in all estimates the value
$e^{-\varepsilon t/2}$
is responsible for the stabilization rate.


\begin{center}
{\bf 5. Conclusion}
\end{center}

In the present paper we have considered a predator--prey model
describing the process of plankton--fish interaction.
We have considered a case of asymptotic stability of the equilibrium point
corresponding to the presence of phytoplankton and zooplankton in the system
and the absence of fish.
The Lyapunov--Krasovskii functional has been constructed,
with the help of which conditions for the initial data have been indicated,
under which plankton survival and fish extinction occur,
and estimates have been obtained,
that characterize the stabilization rate of the amount
of phytoplankton and zooplankton to constant values
and the decrease rate of the number of fish.
The obtained theoretical results can be used for numerical study
of behavior of solutions to the considered model.
\\

The author is grateful to Professor G.V.~Demidenko
and Professor I.I.~Matveeva for the attention to the research.


\end{document}